\documentstyle[amssymb,amsmath,eucal]{article}

\begin{document}


\def\K{{\Bbb K}}
\def\R {{\Bbb R }}
\def\C {{\Bbb C }}
\def\H {{\Bbb H }}

\def\U{{\rm U}}
\def\Sp{{\rm Sp}}
\def\const{{\rm const}}
\def\B{{\rm B}}
\def\O{{\rm O}}
\def\SO{{\rm SO}}
\def\SOS{{\rm SO}^*}
\def\GL{{\rm GL}}
\def\SL{{\rm SL}}
\def\SU{{\rm SU}}

\def\UU{{\cal U}}

\def\Gr{{\rm Gr}}
\def\SS{{\cal S}(B,J)}

\def\ov{\overline}

\def\phi{\varphi}
\def\psi{\varpsi}
\def\epsilon{\varepsilon}
\def\le{\leqslant}
\def\ge{\geqslant}

\def\kvadrat{\hfill{$\blacksquare$}}

\def\sym{ $B$ is symmetric,}
\def\kos{$B$ is skewsym.,}
\def\erm{$B$ is hermitian,}
\def\aerm{$B$ аis antihermitian,}

\def\lin{$J$ is linear,}
\def\alin{$J$ is аantilin.,}

\def\lap{$J^2=1$,}
\def\lam{$J^2=-1$,}
\def\mup{$B(Jv,Jw)=B(v,w)$.}
\def\mum{$B(Jv,Jw)=-B(v,w)$.}
\def\mupa{$B(Jv,Jw)=\ov{B(v,w)}$.}
\def\muma{$ B(Jv,Jw)=-\ov{B(v,w)}$.}

\def\vtch{$\quad$ В In particular\,}
\def\comp{ is a compact symmetric space}
\def\noncomp{ is a noncompact symmetric space}

\def\G{\GL^J}
\def\GG{\UU(D)}

\newcounter{sec}
 \renewcommand{\theequation}{\arabic{sec}.\arabic{equation}}
\newcounter{nomer}

\newcommand{\nomer}{\addtocounter{nomer}{1}%
{\smallskip\noindent$\bf \arabic{nomer}.^{\phantom{\star}}$\ }}
\newcommand{\nomers}{\addtocounter{nomer}{1}%
{\smallskip\noindent$\bf \arabic{nomer}.^{\star}$\ }}

\begin{center}
   {\bf    Pseudoriemannian symmetric spaces: one-type realizations
and open embeddings to grassmannians}

  \vspace{22pt}
 Neretin Yu.A. \footnote{ Supported by RFFI grants 98-01-00303 and
 96-01-96249}

neretin@main.mccme.rssi.ru

\end{center}

Consider a real semisimple group $G$. Let
$\sigma$ be an automorphism of
 $G$ of order 2 (i.e. $\sigma^2=1$). Denote by
$H$ the set of fixed points of the automorphism $\sigma$.
Homogeneous spaces $G/H$
are known as pseudoriemannian {\it symmetric spaces}
(another term  -- affine symmetric spaces; below we name them
by symmetric spaces).
According to Berge classification
 \cite{Ber}
there exists 54 series of classical pseudoriemannian
symmetric spaces and 131 exceptional spaces
\footnote{Of course classification is given
up to coverings and connected components.
A symmetric space $G/H$ is called classical if the group $G$ is
classical (in this case the group $H$ also is classical). 
Number of series in various versions of the list is different.
The reason is preasence of the spaces
 $\O(n+2,\C)/\O(n,\C)\times\O(2,\C)$
(and their real forms). Properties of these spaces
are very specific and hence these spaces often
are not included to the series
$\O(n+m,\C)/\O(n,\C)\times\O(m,\C)$.}.
This list contains many interesting objects from various
branches of mathematics. In particular it contains
riemann symmetric spaces
(including spheres, Lobachevskii spaces, Siegel upper half plane,
grassmannians, matrix balls, Cartan domains, future tubes,
quadrics in $\C{\Bbb P}^n$, symmetric cones, moduli space of
$K3$-surfaces), complex symmetric spaces
(including the space ${\rm PGL}(n,\C) /\O(n,\C)$
of all nondegenerate quadrics), simple Lie groups\footnote%
{Let $P$ be a simple Lie group. Then the group $G=P\times P$
acts on $P$ by left and right multiplications,
the stabilizer of the point $e$ is $H=P$},
multidimensional hyperboloids, spaces of correlations.

Harmonic analysis for many symmetric spaces appeared
long ago. Analysis on sphere and Lobachevskii plane
is a classical subject (from A.M.Legendre to F.G.Mehler, P.Funk,
J.Radon). Analysis on semisimple groups was a subject of investigations
of I.M.Gelfand and M.A.Naimark, Harish-Chandra and their
numerous 
successors.  Analysis on riemann symmetric spaces
was a subject of works of Hua Loo Keng, S.G.Gindikin and F.I.Karpelevich,
S.Helgason. Analysis on hyperboloids was investigated by
V.F.Molchanov. Studying of analysis on arbitrary
pseudoriemannian symmetric spaces was initiated by
M.Flensted-Jensen work about discrete series
 (H.Schlichtkrull, T.Oshima, J.Sekiguchi, V.F.Molchanov, \\ E.P. van den Bahn,
P.Delorm). Now this subject is one of the most wide
and interesting  branches of noncommutative harmonic analysis.

The purpose of this paper is to formulate two simple observations
(by strange way they were not known
\footnote{I think the reason is the length of the Berge list
and nonroot nature of our construction}).
The first observation is:

\smallskip

 {\it all 54 series of classical pseudoriemannian
symmetric spaces have extremely simple uniform geometric realizations.}
Precisely a point of a symmetric space is a pair
of complementary subspaces
 $V_1$, $V_2$ in $\R^{k}$, $\C^{k}$,
$\H^{k}$  satisfying very simple conditions
(isotropy, orthogonality or existence of fixed operator
transposing
$V_1$ и $V_2$).

\smallskip

 It seems to me that this observation pleasantly
simplifies the picture. It also gives a possibility to work with
arbitrary  {\it classical} symmetric spaces by uniform elementary methods.
Several applications of this point of view are in section 
 \S 3.
In particular for all classical symmetric spaces we introduce
simple matrix coordinates (this construction generalizes
Cartan matrix balls
\cite{Car}).

The second observation is: 

{\it all 
classical symmetric spaces have natural open embeddings to grassmannians}. 

This fact was known
for many particular cases. It was many discussed from geometrical
point of view and it was intensively used in harmonic analysis
(see for instance the work of Hua
 \cite{Hua},
more general observations see in papers of
B.O.Makarevich\cite{Mak}, W.Bertram \cite{Bert} and S.G.Gindikin \cite{Gin},  the last paper also contains  some
 bibliography;  I have to apologize to all authors which
are not cited here). We show that this embedding exists
for all%
\footnote{up to covering, center and components of connectedness}
classical symmetric spaces. Moreover in our model it is
absolutely obvious. Direct consequence of this observation
is the following fact (see precise formulation in Section 3):

\smallskip

{\it Representation of the group $G$ in
$L^2(G/H)$ is the restriction of some representation of
some larger group $G^*\supset G$.}

      \medskip

{\bf \S 0. Notations}

      \medskip

Here for avoiding of ambiguity we fix a terminology and notations.

{\bf 0.0.}
The symbol $\K$ denote $\R$, $\C$ or quaternionic ring $\H$.
The term {\it a linear space over $\H$} means (for us) a right module over
$\H$. It is convenient to think that elements of
 $\H^n$ are  vector-columns $v$. We wright linear operators over $\H$
in a form
$v\mapsto Av$
where $A$ is a matrix. We wright a multiplication by a scalar
in a form
$v\mapsto v\lambda $.

 {\it An antilinear operator} in a linear space $V$ over $\C$
is a map $V\mapsto V$ satisfying the conditions
$$A(v+w)=Av+Aw\qquad\qquad A\lambda v=\overline{\lambda}Av,\quad\mbox{где}
\,\,\,\lambda\in\C$$

{\sc  Remark.}  Antilinear maps in linear spaces over $\H$
don't exist. Indeed let $A$ be an antilinear map. Then
\begin{gather*}A(v\lambda\mu)  =A((v\lambda)\mu )=
(A(v\lambda))\ov\mu=(Av)\ov\lambda\ov\mu,\\
 A(v\lambda\mu) = A(v(\lambda\mu))
 =(Av)\ov {\lambda\mu}=(Av)\ov\mu\ov\lambda
\end{gather*}

{\bf 0.1. Forms.}
We remind that a map $B(v,w)$ from $\K^n\times\K^n$ to $\K$ is named
{\it sesquilinear} if for all $v,v_1,v_2, w, w_1,w_2\in\K^n$,
$\lambda\in\K$
we have
\begin{align*}
 B(v\lambda ,w)&= B(v,w)\lambda  &
 B(v_1+v_2, w) &= B(v_1,w)+B(v_2,w)\\
 B( v, w\lambda )&=\overline{\lambda}B(v,w)&
 B(v, w_1+w_2) &= B(v,w_1)+B(v,w_2)
\end{align*}
A sesquilinear map is a {\it hermitian form},
if
$B(v,w)=\overline{B(w,v)}$   для всех $v,w\in\K^n$.
A sesquilinear map is a {\it antihermitian form} if
$B(v,w)=-\overline{B(w,v)}$.

The term {\it form} everywhere in this paper means a  {\it nondegenerate} form
on a linear space over $\K=\R,\C,\H$ having one of the following
7 types:

\smallskip

   --   over $\R$ -- bilinear symmetric and skew symmetric forms

\smallskip

  --   over $\C$ -- bilinear symmetric and antisymmetric forms
and also hermitian and antihermitian forms

\smallskip

  --   over  $\H$  -- hermitian and antihermitian forms

\smallskip

{\sc Remarks.} a) Of course hermitian (antihermitian)
forms over $\R$ are the same as symmetric (skew symmetric)
forms.

b) Antihermitian forms over $\C$ differ unessentialy from
hermitian forms.
Indeed let $B(v,w)$ be a antihermitian form. Then the form $iB(v,w)$
is hermitian.

c)
 Bilinear forms on $\H^n$ don't exist (an expression
$\sum x_sy_s$ is not bilinear form on right $\H$-module!).

d) Hermitian (antihermitian) forms over $\H$ can be represented as
$B(v,w)=\sum \overline w_s b_{st} v_t$ where $b_{st}=\overline b_{ts}$
(resp. $b_{st}=-\overline b_{ts}$).

Remind the classification of forms up to a linear change of variables.
Hermitian forms%
\footnote{including antihermitian forms over $\C$}
over $\R$, $\C$, $\H$ are enumerated by inertia indexes.
In all other cases all nondegenerate forms of the given type
on a given linear space  are equivalent.

{\bf 0.2. Classical groups.} We denote by $\UU(B)$ the group
 of all linear operators
 preserving  a form $B$.
Fix notations for all 7 types of the groups
$\UU(B)$:

 $\Sp(2n,\R)$, $\Sp(2n,\C)$ -- the groups of all linear operators
in $\R^{2n}$,
$\C^{2n}$ preserving skewsymmetric bilinear form.

  $\O(n,\C)$ -- the group of operators in $\C^n$ preserving
 symmetric bilinear form.

$\O(p,q)$, $\U(p,q)$, $\Sp(p,q)$ --
 the groups of operators in  $\R^{p+q}$,
  $\C^{p+q}$, $\H^{p+q}$ preserving a hermitian form with
inertia indexes $(p,q)$

$\SOS(2n)$ -- the group of operators in $\H^n$,
 preserving antihermitian form.

By the  term        {\it "classical group"} we mean a group
of this 7 series and also

 $\GL(n,\R)$, $\GL(n,\C)$,  $\GL(n,\H)$ -- groups of all linear
operators in $V=\R^n, \C^n, \H^n$. We also use notation
$\GL(V)$.

Emphasis that  {\it we don't include to this list groups}
$$\SL(n,\R),\quad \SL(n,\C),\quad \SL(n,\H),\quad \SU(p,q),\quad
 {\rm PSL}(n, \R)\quad , {\rm SL^\pm }(n,\R) \text { etc.}$$
Accordingly we consider (for instance) symmetric spaces
$\U(p,q)/\O(p,q)$ and don't consider  $\SU(p,q)/\SO(p,q)$.
From the point of view of harmonic analysis there is no difference between
spaces
$\U(p,q)/\O(p,q)$ and $\SU(p,q)/\SO(p,q)$.
From the point of view of our realization first space is more
pleasant.

{\bf 0.3. Grassmannians.}  Denote by $\Gr_p(V)$ the set of all
$p$-dimensional subspaces in a linear space $V$.

Consider a form $B$ in $V$.
Remind that a subspace $Q\subset V$ is named
 {\it isotropic} with respect to the form $B$,
if $B$ equals identical zero on $Q$.

Recall that a form $B$ in $V$ is named {\it split}, if there exists
$B$-isotropic half-dimensional subspace. Recall that
a hermitian  form is split if its positive and negative
inertia indexes coincides.
   Orthogonal(symmetric bilinear) forms over
$\C$ and antihermitian forms over $\H$
are split iff dimension of the space is even.
Skew-symmetric forms  always are split.

If
$B$ is split then there exist a basis
$e_1,\dots,e_n,f_1,\dots,f_n$ in the space
$V$  such that
$$B(e_k,e_l)=0,\quad B(f_k,f_l)=0,\quad B(e_k,f_l)=\delta_{k,l}$$

Let $B$ be a {\it split} form.
Denote by $\Gr(V,B)$ a set of all maximal isotropic subspaces
in $V$.

A word {\it grassmannian} below means
$\Gr_p(V)$     or $\Gr(V,B)$.

\medskip

{\bf \S 1.  Semiinvolutions and their centralizer}

\medskip
\addtocounter{sec}{1}
\setcounter{equation}{0}

{\bf 1.1. Semiinvolutions.}
{\it A semiinvolution} (for detailed discussion and
definition for arbitrary division rings see
 Dieudonne book
\cite{Die}) in a linear space $V$ over
 $\K=\R,\C,\H$
is a linear or antilinear operator $J$ satisfying the condition
 $$J^2=\lambda$$
where $\lambda$ is an element of center of $\K$.

{\sc Lemma 1.1.} {\it Let $J$ be an antilinear
semiinvolution over $\C$,
let $J^2=\lambda$.
Then $\lambda\in\R$.}

{\sc Proof.} Calculate $J^3$ by two ways:
$$J^3v=JJ^2v=J\lambda v=\ov\lambda v\qquad
J^3v=J^2Jv=\lambda Jv \qquad\qquad\qquad\qquad\qquad\qquad
\shoveright{\blacksquare}
$$

Consider a map $S(J):\Gr_p(V)\to \Gr_p(V) $ given by the formula
$Q\mapsto JQ$. Obviously $S(J)^2=1$. The main object of our
interest are these maps.
For each element $\sigma$
of center of $\K$ we have $S(\sigma J)=S(J)$. By this
reason we assume
$$\boxed{J^2=\pm1}$$
Moreover if $\K=\C$ and  $J$ is linear then we can assume $J^2=1$.

Denote by      $\GL^J=\GL^J(V)$  the centralizer of the semiinvolution
$J$  in the group
$\GL(V)$.

{\bf 1.2. Description of semiinvolutions.} A semiinvolution
$J$ defines an additional structure in a linear space
 $V$. This structure is discussed below.

\smallskip

a) Let $J$ be linear ,and $J^2=1$. Consider the subspaces
$V_\pm\subset V$ which consist of vectors $v$ satisfying conditions
$Jv=\pm v$. Then $V=V_+\oplus V_-$. Thus the semiinvolution $J$
defines the  fixed decomposition of $V$ into a direct
sum of two subspaces. Obviously
$\GL^J=\GL(V_+)\times\GL(V_-)$.

\smallskip

b) Let $\K=\R$, $J^2=-1$, $V=\R^{2n}$. Then we can consider the
space $V$ as a space over field $\C$ where multiplication
by a scalar $i$ is the operator
$J$.  Obviously $\GL^J\simeq\GL(n,\C)$.

\smallskip

c) Let $\K=\C$, $J$ is antilinear, $J^2=-1$.
Define the action of the algebra $\H$ on $V$. For these purpose
we assume that the subalgebra $\C\subset \H$ acts as it acts
and that the quaternionic imaginary unit
$j$ is the operator $J$.

Thus   we define a structure of
a linear space over $\H$ on  the space $V$. Clearly   the group
$\GL^J$ is a quaternionic group $\GL$.

\smallskip

d) Let $\K=\C$, $\dim V=k$, let $J$ be antilinear, $J^2=1$.
Consider the set $W$ of all fixed points of the involution
$J$. Obviously the set $W$ is a $\R$-subspace in
$V$ and $V=V\oplus iV$.  Thus we can consider $V$ as a complexification
of the real linear space
$W$.
Obviously $\GL^J\simeq\GL(k,\R)$.

\smallskip

e) Let $\K=\H$, $\dim V=k$, $J^2=-1$.
 Consider the set $W$ of all points $v\in V$ satisfying the condition
$$  Jv=vi$$
It is clear that $W$ is a linear space over $\C$ and
$V=W\oplus Wj$. Hence the space $V$ is a "quaternionization"
of $k$-dimensional complex space $W$ (i.e.
$V=W\otimes_\H V$). Obviously $\GL^J\simeq\GL(k,\C)$.

\smallskip

We say that semiinvolution is {\it split} if there exists a subspace
$Q$ such that $V=Q\oplus JQ$. Only this case is ineresting for us.
 All semiinvolutions
of types b)--e)
always are split. An involution of the type a) is split iff
$\dim V_+=\dim V_-$.

{\bf 1.3. Semiinvolutions consistent with forms.}
Let $B$ be a form, let $J$ be a {\it linear} semiinvolution.
We say that the semiinvolution
$J$
is consistent with the form
$B$ if for all $v,w$ we have
\begin{equation}B(Jv,Jw)=\mu B(v,w)\end{equation}
where $\mu$ is an element of the center of the division ring $\K$.

We say that an {\it antilinear} semiinvolution $J$ is consistent with
a form
$B$ if the following condition   is fulfilled
\begin{equation}B(Jv,Jw)=\mu\ov{ B(v,w)}\end{equation}
Now we will show that $\mu=\pm1$

{\sc Lemma 1.2.} a) {\it Let $J$ be linear. Then $\mu=\pm1$.}

b){\it Let $J$ be antilinear and the form $B$ is hermitian.
Then $\mu=\pm1$.}

c){\it Let $J$ be antilinear and the form $B$ is bilinear. Then
there exists a constant $\sigma$ such that  for the form $\sigma B$
we have $\mu=1$.}

{\sc Proof.} a) Bearing in mind the condition $J^2=\pm1$ we obtain
$$B(v,w)=B(J^2v,J^2w)=\mu B(Jv, Jw)=\mu^2 B(v,w)$$

b) $B(v,w)=B(J^2v,J^2w)=\mu \ov{B(Jv, Jw)}=\mu\ov\mu B(v,w)$

\noindent Hence $|\mu|=1$. Further we substitute $v=w$ to (1.2) and observe
that $B(v,v)$ is real.

c) Similarly we obtain $|\mu|=1$.
Assume $C(v,w)=\frac1\sigma B(v,w)$. Then
$$\sigma C(v,w)=\mu\ov\sigma \ov{ C(v,w)}$$
and now the statement becomes obvious. \kvadrat

{\sc Lemma 1.3.} {\it Let $J$ be a semiinvolution consistent with
the form
 $B$. Let $Q$ be a subspace isotropic with respect
to the form
$B$. Then the subspace $JQ$ is isotropic with respect to $B$.}

{\sc Proof.} If $B(v,w)=0$ then $B(Jv,Jw)=0$. \kvadrat

\smallskip

Thus {\it the semiinvolution $J$ defines an involution on
$B$-isotropic grassmannian.}

{\bf 1.4. Managing form.} Fix a form $B$ and a semiinvolution which is
consistent   with
 $J$. Consider the expression
$$D(v,w):=B(v,Jw)$$
It is easy to check (in each particular case it is obvious) that
$D(v,w)$  is a "form" in our sense.

{\sc Example.} Let $B$ be a symmetric bilinear form over $\R$.
Fix $\lambda$ and  $\mu$. Then
$$D(v,w)=B(v,Jw)=\mu B(Jv, J^2 w)=\mu\lambda B(Jv,w)=\mu\lambda B(w,Jv)=
\mu\lambda  D(w,v)$$
and the type of the form $D$ becomes obvious.

{\sc Lemma 1.4.} {\it Let a subspace  $Q$ be isotropic with respect
to the form
$B$. Then the subspace 
$JQ$ is orthogonal $Q$ with respect to managing form $D$.
}

{\sc Proof.} Let $v,w\in Q$. Then $D(v,Jw)=B(v,J^2w)=0$.\kvadrat

Hence we have 3 structures in the space $V$: two forms and the
 semiinvolution.
If we know two of these structures we can reconstruct the third structure.
We can formulate this remark in the following form:

Denote by
$\UU(B)^J$ the set of all elements of $\UU(B)$ commuting with $J$. Then
\begin{align}
  \UU(B)^J& =\UU(B) \cap  \GL(V)^J\\
  \UU(B)^J& =\UU(B)     \cap   \UU(D) \\
  \UU(B)^J& = \UU(D)    \cap  \GL(V)^J
\end{align}

{\bf 1.5. Split pairs $(B,J)$.} Let $J$  be a semiinvolution
which is consistent with a form
 $B$. We say that the pair $(B,J)$ is {\it split}, if
there exists an isotropic subspace $Q$ in $V$ such that
 $V=Q\oplus JQ$.  We also say that $V$ is a {\it spiting} subspace.

{\sc Lemma 1.5.} {\it Let $(B,J)$ be a split pair. Then for each
 maximal isotropic subspace $P$ in $V$
the subspace $JP$ coincides with  $D$-orthogonal complement to $P$.}

{\sc Proof.} The statement follows from Lemma 1.4  and
calculation of dimensions.    \kvadrat

{\sc Proposition 1.6.}{\it Fix a type of a form $B$, the type of
a semiinvolution $J$
and number $\mu=\pm 1$.  Then there exists an unique
{\rm(} up to linear change of coordinates{\rm )}
a split pair $(B,J)$ of the given type.}

{\sc Proof.} It is easy to see that the pair $(B,J)$
is uniquely defined by the restriction of
managing form
$D$ to the subspace $Q$. Now the problem is reduced to
classification of forms.
 \kvadrat

{\bf 1.6. Description of groups $\UU^J(B)$ for split pairs $(B,J)$.}
In this subsection we give long enumeration of species parallel
to enumeration of Subsection 1.2. We preserve numeration a),b),c)
etc. from Subsection 1.2.

In fact it is possible to carry out the same work in general form
for arbitrary division rings and not only for split pairs
$(B,J)$,
see. (\cite{Die}).

\smallskip

a) Let $J$ be linear, $J^2=1$. Obviously subspaces
$V_\pm$ are $\UU^J(B)$-invariant. Hence
elements of the group $\UU^J(B)$ have block structure
\begin{equation}
\begin{pmatrix} g_1&0\\0&g_2\end{pmatrix}
\end{equation}
We consider two subcases.

\smallskip

a1) Let $\mu=+1$. Then subspaces $V_+$,  $V_-$ are orthogonal
with respect to the form
$B$. Denote by $B_\pm$ the restriction of the form
 $B$ to $V_\pm$.
Obviously   $\UU^J(B)=\UU(B_+)\times \UU(B_-)$.

 A spliting subspace $P$ is a graph of invertible operator
$L_P:V_+\to V_-$. Clearly
operator $L_P$ identifies the form $B_+$ with the form $(-B_-)$
(it is equivalent to isotropy of subspace $P$). Hence
$\UU(B_+)\simeq \UU(B_-)$.

\smallskip

a2) Assume $\mu=-1$. Let $v_1,v_2\in V_+$. Then
$$B(v_1,v_2)=-B(Jv_1,Jv_2)=-B(v_1,v_2)$$
Hence the subspaces $V_\pm$ are $B$-isotropic. The form
$B$ defines nondegenerate pairing between $V_+$ and $V_-$.
Hence the operator
$g_1$  (see  (1.6)) is contragredient  to $g_2$
with respect to our pairing.

Thus $\UU^J(B)=\GL(V_+)$

\smallskip

b) Let $\K=\R$, $J^2=-1$. Then our space over $\R$
can be considered as a space over $\C$. Let us define in this space
a $\C$-valued form
$$Z(v,w)=B(v,w)+iD(v,w)=B(v,w)+iB(v,Jw)$$
Now $\UU^J(B)=\UU(Z)$.

\smallskip

c) Let $\K=\C$, let  $J$ be antilinear, $J^2=-1$.
Let us define in our $\H$-space
a $\H$-valued form
$$Y(v,w)=B(v,w)+jD(v,w)=B(v,w)+jB(v,Jw)$$
We have $\UU^J(B)=\UU(Y)$.

\smallskip

d) Let $\K=\C$, let $J$ be antilinear, $J^2=1$.
  Consider the $\R$-subspace
$V_\R$ consisting of fixed points of the semiinvolution $J$.
Consider the restriction of the form $B$ to $V_\R$. For all
$v,w\in V_\R$ we have
$$B(v,w)=\mu \ov{B(Jv,Jw)}=\mu\ov{ B(v,w)}$$
Now define the form $X$ on $V_\R$ which equals $B$ if $\mu=1$,
and equals $iB$ if $\mu=-1$. Then the form $X$ is a
$\R$-value form and $\UU^J(B)=\UU(X)$.

\smallskip

e) Let $\K=\H$, $J^2=-1$. This case is similar to d).

\medskip
\begin{center}
        {\bf \S 2. Realizations of classical symmetric spaces}

\medskip
\end{center}

\addtocounter{sec}{1}
\setcounter{equation}{0}

In fact models of all classical symmetric spaces were obtained
in previous section. We have only to write the list.

              {\bf 2.1. List 1. The case of
split pairs $(B,J)$.}
Fix a space $V=\K^\alpha$ and split pair $(B,J)$ in this space (see 1.5).
We will say that the form $B$ is an  {\it underlying form},
and a semiinvolution  $J$ is a managing {\it
semiinvolution}.
Let $Q$ be the  associated managing form, see 1.4.

Let us define the space $\SS$. Its points are ordered pairs
of subspaces
$Q_1$, $Q_2$ in $V$ satisfying the conditions

$1$. $Q_1$, $Q_2$ are maximal $B$-isotropic subspaces

$2$. $V=Q_1\oplus Q_2$.

$3$. $JQ_1=Q_2$ or equivalently $JQ_2=Q_1$ (or
equivalently $Q_2$ is $D$-orthogonal complement to $Q_1$)

Consider the centralizer $G=\UU^J(B)$ of the semiinvolution
$J$ in the group $\UU(B)$. The group $G=\UU^J(B)$ acts
on the space $\SS$ by the obvious way.
We claim that either
$\SS$ is a symmetric space $G/H$ or
$\SS$ is an union of finite family of symmetric spaces
of the type $G/H_i$.

The group $G$ was described in Subsection 1.6.
Now we want to describe the stabilizer $H$
of the pair of subspaces $(Q_1,Q_2)$. For this purpose we define
the form
$D'(v,w)=B(v,Jw)$ on subspace
$Q_1$ (it is the restriction of managing form $D$ to
the subspace $Q_1$). It is easy to see that
$H\simeq\UU(D')$.

Below we give a list of symmetric spaces obtained in this way.

\smallskip

The first line indicate the symmetric space
 $G/H$.

\smallskip

The second line contains the space $V=\K^l$ and
 type of the form $B$. In this line we also indicate 
 type of the semiinvolution $J$

\smallskip

The third line contains the group $\UU(B)$. For us it is more
pleasant to denote it by $G^*$
(in some cases this information also give precise inertia
indexes
of the form $B$).

\smallskip

The forth line contains the centralizer
 $\G$ of the semiinvolution $J$ in $\GL(V)$ and also the group
$\GG$ consisting of operators preserving managing form
 $D(v,w)=B(v,Jw)$.  Call to mind that equalities (1.3)%
--(1.5) are fulfilled.

\smallskip

We mark by the symbol $\star$ the cases when
$\SS$  is not $G$-homogeneous space.
In this case we give a decomposition of
 $\SS$ onto an union of symmetric spaces.

\pagebreak

\def\table{\begin{tabbing}
$\quad V=R^{2p(+q)}$,\= $B$ kososimmetr.,\=
        $J$  antilinejna,\= $J^2=-1$, \=\muma\kill}
\def\etable{\end{tabbing}}

\table
\nomer $\O(p,q)\times\O(p,q)/\O(p,q)$\\
$\quad V=\R^{2(p+q)}$,\> \sym\> \> \lap\> \mup\\
$\quad G^*=\O(p+q,p+q)$,\\
 $\quad\G=\GL(p+q,\R)\times\GL(p+q,\R)$\>\>\>
                                     $ \GG=\O(2p,2q)$.\\
\vtch $\O(p)\times \O(p)/\O(p)$ \comp
\etable

\table
\nomer $\Sp(2n,\R)\times \Sp(2n,\R)/\Sp(2n,\R)$\\
$\quad V=\R^{4n}$,\> \kos\> \> \lap\> \mup\\
$\quad G^*=\Sp(4n,\R)$                            \\
$\quad \G= \GL(2n,\R) \times \GL(2n,\R)$,\>\>\> $\GG=\Sp(4n,\R)$
\etable

\table
\nomers $\GL(n,\R)/\O(p,n-p)$\\
$\quad V=\R^{2n}$,\>\kos\>\>\lap\>\mum \\
$\quad G^*=\Sp(2n,\R)$\\
$\quad \G=\GL(n,\R)\times\GL(n,\R)$\>\>\>$\GG=\O(n,n)$\\
$\quad \SS=\cup_{p=0}^n  \GL(n,\R)/\O(p,n-p) $\\
 \vtch $\GL(n,\R)/\O(n)$ \noncomp
\etable

\table
\nomer $\GL(2n,\R)/\Sp(2n,\R)$\\
$\quad V=\R^{4n}$,\>\sym\>\> \lap\>\mum\\
$\quad G^*=\O(2n,2n)$ \\
$\quad \G=\GL(2n,\R)\times \GL(2n,\R)$, \>\>\>$\GG=\Sp(2n,\R)$
\etable

\table
\nomers $\O(n,\C)/\O(p,n-p)$              \\
$\quad V=\R^{2n}$, \>\sym \>\> \lam\> \mum\\
$\quad G^*=\O(n,n)$\\
$\quad \G=\GL(n,\C)$,\>\>\>$\GG=\O(n,n)$\\
$\quad \SS=\cup_{p=0}^n \O(n,\C)/\O(p,n-p) $\\
 \vtch $\O(n,\C)/\O(n)$ \noncomp
\etable

\table
\nomer $\Sp(2n,\C)/\Sp(2n,\R)$\\
$\quad V=\R^{4n}$,\>\kos\>\>\lam\>\mum\\
$\quad G^*=\Sp(4n,\R)$\\
$\quad \G=\GL(2n,\C)$\>\>\>$\GG=\Sp(4n,\R)$
\etable

\table
\nomer $\U(n,n)/\Sp(2n,\R)$\\
$\quad V=\R^{4n}$, \>\sym\>\>\lam\>\mup\\
$\quad G^*=\O(2n,2n)$\\
$\quad \G=\GL(2n,\C)$\>\>\>$\GG=\Sp(4n,\R)$
\etable

\table
\nomer $\U(p,q)/\O(p,q)$\\
$\quad V=\R^{2(p+q)}$,\>\kos\>\>\lam\>\mup\\
$\quad G^*=\Sp(2(p+q),\R)$\\
$\quad \G=\GL(p+q,\C)$\>\>\> $\GG=\O(2p,2q)$\\
 \vtch $\U(p)/\O(p)$ \comp
\etable

\table
\nomer $\O(n,\C)\times \O(n,\C)/\O(n,\C) $ \\
$\quad V=\C^{2n}$,\>\sym\>\lin\>\lap\>\mup\\
$\quad G^*=\O(2n,\C)$,\\
$\quad \G=\GL(n,\C)\times\GL(n,\C)$\>\>\>$\GG=\O(2n,\C)$
\etable

\table
\nomer $\Sp(2n,\C)\times \Sp(2n,\C)/\Sp(2n,\C)$\\
$\quad V=\C^{4n}$,\> \kos \> \lin\> \lap\>\mup \\
$\quad G^*=\Sp(4n,\C)$\\
$\quad \G=\GL(2n,\C)\times\GL(2n,\C)$\>\>\>$\GG=\Sp(4n,\C)$
\etable

\table
\nomer $\U(p,q)\times \U(p,q)/ \U(p,q)$\\
$\quad V=\C^{2(p+q)}$,\> \erm \> \lin \>\lap \> \mup\\
$\quad G^*=\U(p+q,p+q)$\\
$\quad \G=\GL(p+q,\C)\times\GL(p+q,\C)$\>\>\>$\GG=\U(2p,2q)$\\
\vtch $\U(p)\times\U(p)/\U(p)$ \comp
\etable

\table
\nomer $\GL(2n,\C)/\Sp(2n,\C)$\\
$\quad V=\C^{4n}$, \> \sym \> \lin \> \lap \>  \mum \\
$\quad G^*=\O(4n,\C)$,\\
$\quad \G=\GL(2n,\C)\times \GL(2n,\C)$\>\>\>$\GG=\Sp(4n,\C)$
\etable

\table
\nomer $\GL(n,\C)/\O(n,\C)$\\
$\quad V=\C^{2n}$, \> \kos\> \lin \> \lap \>\mum\\
$\quad G^*=\Sp(2n,\C)$\\
$\quad \G=\GL(n,\C)\times\GL(n,\C)$\>\>\>$\GG=\O(2n,\C)$.
\etable

\table
\nomers $\GL(n,\C)/\U(p,n-p)$ \\
$\quad V=\C^{2n}$, \> \erm \> \lin\> \lap \> \mum\\
$\quad G^*= \U(n,n)$\\
$\quad \G=\GL(n,\C)\times\GL(n,\C)$,\>\>\>$\GG=\U(n,n)$\\
$\quad \SS=\cup_{p=0}^n \GL(n,\C)/\U(p,n-p)$\\
\vtch $\GL(n,\C)/\U(n)$ \noncomp
\etable

\table
\nomers $\Sp(2n,\R)/\U(p,n-p)$\\
$\quad V=\C^{2n}$,\>\kos \>\alin \> \lap\> \mupa\\
$\quad G^*=\Sp(2n,\C)$,\\
$\quad \G=\GL(2n,\R)$\>\>\> $\GG=\U(n,n)$\\
$\quad \SS=\cup_{p=0}^n\Sp(2n,\R)/\U(p,n-p)$\\
\vtch $\Sp(2n,\R)/\U(n)$ \noncomp
\etable

\table
\nomer $\O(2p,2q)/\U(p,q)$\\
$\quad V=\C^{2(p+q)}$,\> \sym \> \alin \> \lap\> \mupa \\
$\quad G^*=\O(2(p+q),\C)$\\
$\quad \G=\GL(2(p+q),\R)$\>\>\>$\GG=\U(2p,2q)$\\
\vtch $\O(2p)/\U(p)$ \comp
\etable

\table
\nomer $\O(n,n)/\O(n,\C)$ \\
$\quad V=\C^{2n}$, \> \erm\> \alin \> \lap\> \mupa \\
$\quad G^*=\U(n,n)$, \\
$\quad \G= \GL(2n,\R)$\>\>\> $\GG= \O(2n,\C)$
\etable

\table
\nomer $\Sp(4n,\R)/\Sp(2n,\C)$ \\
$\quad V=\C^{4n}$, \> \erm \> \alin \> \lap\> \muma \\
$\quad G^*=\U(2n,2n)$,\\
$\quad\G=\GL(4n,\R)$,\>\>\> $\GG=\Sp(4n,\C)$
\etable

\table
\nomer $\Sp(p,q)/\U(p,q)$  \\
$\quad V=\C^{2(p+q)}$, \> \kos \> \alin\> \lam \> \mupa\\
$\quad G^*=\Sp(2(p+q),\C)$\\
$\quad \G=\GL(p+q,\H)$, \>\>\> $\GG=\U(2p,2q)$\\
\vtch $\Sp(p)/\U(p)$ \comp
\etable

\table
\nomers $\SOS(2n)/\U(p,n-p)$\\
$\quad V=\C^{2n}$,\> \sym\> \alin \> \lam\> \mupa\\
$\quad G^*=\O(2n,\C)$,\\
$\quad \G=\GL(n,\H)$\>\>\> $\GG=\U(n,n)$\\
$\quad \SS=\cup_{p=1}^n \SOS(2n)/\U(p,n-p)$\\
\vtch $\SOS(2n)/\U(n)$ \noncomp
\etable

\table
\nomer $\Sp(n,n)/\Sp(2n,\C)$\\
$\quad V=\C^{4n}$, \>\erm \> \alin \> \lam \> \mupa    \\
$\quad G^*=\U(2n,2n)$\\
$\quad \G=\GL(2n,\H)$\>\>\> $\GG= \Sp(4n,\C)$
\etable

\table
\nomer $\SOS(2n)/\O(n,\C)$\\
$\quad V=\C^{2n}$,\> \erm\> \alin \> \lam \> \muma\\
$\quad G^*=\U(n,n)$\\
$\quad G=\GL(n,\H)$, \>\>\>$\GG= \O(2n,\C)$
\etable

\table
\nomer $\Sp(p,q) \times \Sp(p,q) /\Sp(p,q)  $\\
$\quad V=\H^{2(p+q)}$,\> \erm \>\> \lap\> \mup\\
$\quad G^*=\Sp(p+q,p+q)$, \\
$\quad\G=\GL(p+q,\H)\times\GL(p+q,\H) $\>\>\> $\GG=\Sp(2p,2q)$\\
\vtch $\Sp(p)\times \Sp(p)/\Sp(p)$ \comp
\etable

\table
\nomer $\SOS(2n)\times \SOS(2n)/\SOS(2n) $   \\
$\quad V=\H^{2n}$, \> \aerm \>\> \lap\> \mup\\
$\quad G^*=\SOS(4n)$\\
$\quad \G=\GL(n,\H)\times \GL(n,\H)$\>\>\>$\GG=\SOS(4n)$
\etable

\table
\nomer $\GL(n,\H)/\SOS(2n)$\\
$\quad V=\H^{2n}$,\> \erm \> \> \lap \> \mum\\
$\quad G^*=\Sp(n,n)$,\\
$\quad \G= \GL(n,\H) \times \GL(n,\H)$, \>\>\> $\GG=\SOS(4n)$
\etable

\table
\nomers $\GL(n,\H)/\Sp(p,n-p)$\\
$\quad V=\H^{2n}$, \> \aerm \>\> \lap \> \mum\\
$\quad G^*=\SOS(4n)$\\
$\quad \G=\GL(n,\H) \times \GL(n,\H) $\>\>\> $\GG=\Sp(n,n)$\\
$\quad \SS=\cup_{p=0}^n  \GL(n,\H)/\Sp(p,n-p)$\\
\vtch $\GL(n,\H)/\Sp(n)$ \noncomp
\etable

\table
\nomer $\U(2p,2q)/\Sp(p,q)$\\
$\quad V=\H^{2(p+q)}$,\> \aerm  \>\>  \lam \> \mup\\
$\quad G^*=\SOS(4(p+q))$, \\
$\quad \G= \GL(2(p+q),\C)$,\>\>\> $\GG=\Sp(2p,2q)$\\
\vtch $\U(2p)/\Sp(p)$ \comp
\etable

\table
\nomer $\U(n,n)/\SOS(2n)$\\
$\quad V=\H^{2n}$, \> \erm \>\> \lam\> \mup\\
$\quad G^*=\Sp(n,n)$,\\
$\quad \G=\GL(2n,\C) $,\>\>\>$\GG=\SOS(4n) $
\etable

\table
\nomer $\O(2n,\C)/\SOS(2n)$\\
$\quad V=\H^{2n}$,\> \aerm \>\> \lam \> \mum \\
$\quad G^*= \SOS(4n)$\\
$\quad \G=\GL(2n,\C)$\>\>\>$\GG=\SOS(4n)$
\etable

\table
\nomers $\Sp(2n,\C)/\Sp(p,n-p)$\\
 $\quad V=\H^{2n}$, \> \erm \>\> \lam \> \mum \\
$\quad G^*=\Sp(n,n)$\\
$\quad \G=\GL(2n,\C)$\>\>\> $\GG= \Sp(n,n)$\\
$\quad \SS=\cup_{p=0}^n    \Sp(2n,\C)/\Sp(p,n-p)  $\\
\vtch $\Sp(2n,\C)/\U(n)$ \noncomp
\etable

{\bf 2.2. List 2. The case when we have only
underlying form $B$.}
Let us fix a linear space $V=\K^{2n}$ equipped with a split
form $B$.
Define the space ${\cal S}(B)$. Points of ${\cal S}(B)$
are ordered pairs $(Q_1,Q_2)$ of maximal isotropic subspaces
in the $V$ such that $V=Q_1\oplus Q_2$. Obviously all spaces
${\cal S}(B)$ are symmetric spaces having the type
 $$G/H=\UU(B)/\GL(n,\K)$$

Below we give the list of symmetric spaces obtained in this way

\nomer$\O(n,n)/\GL(n,\R)$

\nomer$\Sp(2n,\R)/\GL(n,\R)$

\nomer$\O(2n,\C)/\GL(n,\C)$

\nomer$\Sp(2n,\C)/\GL(n,\C)$

\nomer$\U(n,n)/\GL(n,\C)$

\nomer$\Sp(n,n)/\GL(n,\H)$

\nomer$\SOS(4n)/ \GL(n,\H)$

In all cases we define the group
$$G^*=G\times G$$

{\bf 2.3. List 3. The case when we have only a managing
semiinvolution.} Consider a linear space $V=\K^{2n}$
and a split (see 1.2) semiinvolution $J$ in $V$.
 Let us define the space ${\cal S}(J)$. Points of ${\cal S}(J)$
are all pairs of subspaces $(Q_1,Q_2)$ in  $V$ such that
$V=Q_1\oplus Q_2$, $JQ_1=Q_2$.
   The group $G:=\GL^J$(centralizer of the semiinvolution $J$)
acts on ${\cal S}(J)$. It is readily seen that
in all cases
 ${\cal S}(J)$
is a symmetric space.

The list of such spaces is given below.
The first row indicates the space $G/H$.
The second row contains the space $V$ and the type of
the semiinvolution $J$.

\def\tabl{\begin{tabbing}
$\quad V=\H^{2n}\qquad$\= $J^2=-1$ ff,\= $J$ linejna\kill }
\def\etabl{\end{tabbing}}

\tabl
\nomer$\GL(n,\R)\times \GL(n,\R) /\GL(n,\R)$\\
$\qquad V=\R^{2n}$, \> $J^2=1$\\

\nomer$\GL(n,\C)/\GL(n,\R)$\\
$\qquad V=\R^{2n}$, \> $J^2=-1$ \\

\nomer$\GL(n,\C)\times\GL(n,\C) / \GL(n,\C) $\\
$\qquad V=\C^{2n}$, \> $J^2=1$, \> $J$ линейна\\

\nomer$\GL(2n,\R)/\GL(n,\C)$\\
$\qquad V=\C^{2n}$\> $J^2=1$, \> $J$ антилинейна\\

\nomer$\GL(n,\H)/\GL(n,\C)$\\
$\qquad V=\C^{2n}$, \> $J^2=-1$,\>  $J$ антилинейна\\

\nomer$\GL(n,\H) \times \GL(n,\H)/\GL(n,\H) $\\
$\qquad V=\H^{2n}$,\> $J^2=1$\\

\nomer$\GL(2n,\C)/\GL(n,\H)$\\
$\qquad V=\H^{2n}$, $J^2=-1$
\etable

Во всех случаях определим группу
$$G^*=\GL(V)$$

{\bf 2.4. List 4. The case when we have  only
managing form.} Consider a space $V=\K^n$,
equipped with a form $D$. Assume $G=\UU(D)$.
 Consider the space ${\cal S}_m(D)$ consisting of
all $m$-dimensional subspaces  $Q_1\subset V$
such that the form $D$ is nondegenerate on $Q_1$.
We also can say that a point of the space ${\cal S}_m(D)$
is a pair of subspaces $(Q_1,Q_2)$ such that

1.$V=Q_1\oplus Q_2$,

2. $Q_2$ is $D$-orthogonal complement to $Q_1$.

Obviously, either ${\cal S}_m(D)$ is a symmetric space
of the type $G/H$
or ${\cal S}_m(D)$ is an union of a finite family of symmetric spaces
$G/H_i$. The list of such symmetric spaces is
given below.
In the case when the form $D$ is hermitian we also describe
decomposition of ${\cal S}_m(D)$ onto the union of symmetric spaces.

\smallskip
\noindent\nomers     $\O(p,q)/\O(r,s)\times\O(p-r,q-s)$

 $\qquad{\cal S}_m(D)=\bigcup\limits_{r,s:\, r+s=m, r\le p, s\le q}
                     \O(p,q)/\O(r,s)\times\O(p-r,q-s)$

 В In particular $\O(p,q)/\O(p)\times \O(q)$ \noncomp

\hspace{1.5cm}     $\O(p)/\O(r)\times\O(p-r)$ \comp

\noindent\nomer $\Sp(2(k+l),\R)/\Sp(2k,\R)\times\Sp(2l,\R)$

\noindent\nomer $\O(n+m,\C)/\O(n,\C)\times\O(m,\C)$

\noindent\nomer $\Sp(2(k+l),\C)/\Sp(2k,\C)\times\Sp(2l,\C)$

\noindent\nomers     $\U(p,q)/\U(r,s)\times\U(p-r,q-s)$\\

 $\qquad{\cal S}_m(D)=\bigcup\limits_{r,s:\, r+s=m, r\le p, s\le q}
                     \U(p,q)/\U(r,s)\times\U(p-r,q-s)$

  В In particular $\U(p,q)/\U(p)\times\U(q)$ \noncomp

 \hspace{1.5cm} $\U(p)/\U(r)\times\U(p-r)$ \comp

\noindent\nomer     $\Sp(p,q)/\Sp(r,s)\times\Sp(p-r,q-s)$\\

 $\qquad{\cal S}_m(D)=\bigcup\limits_{r,s:\, r+s=m, r\le p, s\le q}
                     \Sp(p,q)/\Sp(r,s)\times\Sp(p-r,q-s)$

  В In particular $\Sp(p,q)/\Sp(p)\times\Sp(q)$ \noncomp

 \hspace{1.5cm} $\Sp(p)/\Sp(r)\times\Sp(p-r)$ \comp

\noindent\nomer    $\SOS(2(m+n))/\SOS(2m)\times\SOS(2n)$

In all cases we define the group
$$G^*=\GL(V)$$

{\bf 2.5. List 5. The case when there is nothing.}
Consider the space $V=\K^{p+q}$.
 Further consider the space ${\cal S}_p$,
consisting of all pairs of subspaces $(Q_1,Q_2)$ in $V$
 such that

1. $\dim Q_1=p$, $\dim Q_2=q$

2. $V=Q_1\oplus Q_2$

By this way we obtain the following symmetric spaces $\cal S$

\nomer $\GL(p+q,\R)/\GL(p,\R)\times\GL(q,\R)$

\nomer $\GL(p+q,\C)/\GL(p,\C)\times\GL(q,\C)$

\nomer $\GL(p+q,\H)/\GL(p,\H)\times\GL(q,\H)$

In all cases we define the group
 $$G^*=\GL(p+q,\K)$$

\medskip

{\bf \S 3. Some applications}

\medskip

\addtocounter{sec}{1}
\setcounter{equation}{0}

{\bf 3.1. Open embeddings to grassmannians.}
Thus in all 54 cases a point of a symmetric space $G/H$
 is a pair of subspaces $(Q_1,Q_2)$ in a linear space.

If we have managing semiinvolution or managing form
(Lists 1,3,4),
then the subspace $Q_2$ is uniquely defined by the subspace $Q_1$.
Hence the map $(Q_1,Q_2)\mapsto Q_1$ is an open embedding
of the symmetric space ${\cal S}=G/H$ to some grassmannian $\Gr^*$
(this grassmannian is complete grassmannian for the Lists
3,4 and isotropic grassmannian $\Gr(D)$ for the List
1).

Let we have no  managing semiinvolution and no managing form (Lists 2,5).
 Then
$(Q_1,Q_2)$ is a point of products of two grassmannians.
We also denote this product of grassmanianns by
 $\Gr^*$.

The image of the space $G/H$ in the grassmannian (or product of
two grassmannians) $\Gr^*$
 in all cases is open.
Moreover the image is open in all cases except 10 series marked by the
symbol $\star$.

{\sc Remark.} If a space $G/H$ is compact then
its image coincides with grassmannian.
In other words we realized all 10 series of compact symmetric spaces
as grassmannians.

{\bf 3.2. Overgroup.} For all symmetric spaces $G/H$ we indicated
the group $G^*\supset G$. By the construction the group $G^*$ acts
transitively on the grassmannian $\Gr^*$,
A stabilizer of a point is a maximal parabolic subgroup
in $G^*$.

{\bf 3.3. Restriction from degenerated principal series.}Рассмотрим
Consider the natural unitary
representation $\rho$ of the group $G^*$ in the space $L^2$ on $\Gr^*$.

{\sc Proposition 3.1.} a) {\it For all classical symmetric spaces
except the cases $G/H$  marked by $\star$
the restriction of the representation $\rho$ to the subgroup $G$
is equivalent to the representation of $G$ in $L^2(G/H)$.}

\smallskip 

  b) {\it For cases marked by $\star$ the restriction of $\rho$ to $G$
is equivalent to the representation of $G$ in $\bigoplus L^2(G/H_i)$ (where
the spaces $G/H_i$ are indicated in List.}

{\sc Proof.}  It is an obvious consequence from Subsection 3.1. \kvadrat

{\sc Remark.} Consider the case $G^*=G\times G$ ( Lists 2,5).
Consider the representation of the group $G^*=G\times G$
 in the space $L^2$ on the product of two grassmannians.
Obviously this representation is a tensor product of
two representations of the group $G$.
Hence in this cases the representation of
$G$ in $L^2(G/H)$ is a tensor product of two representations
of the group
 $G$ of degenerated principal series.

{\bf 3.4. Matrix coordinates on symmetric spaces.}
Consider a linear space $V$ and a pair of subspaces
$X$, $Y$ such that $Z=X\oplus Y$. Let $\dim X=\alpha$.
As before we denote by $\Gr_\alpha(V)$ the grassmannian
of all $\alpha$-dimensional subspaces.
Assume $R\in \Gr_\alpha(Z)$ doesn't intersect withс $Y$.
Then $R$ is a graph
of some operator $X\to Y$. This operator is named
{\it an angular operator} of subspace $R$,
 associated with the decomposition $V=X\oplus Y$.

Fix a decomposition $V=X\oplus Y$. We wrighte elements of the group
$\GL(V)$  as block operators
$$
\begin{pmatrix} A&B\\C&D\end{pmatrix}:X\oplus Y  \to X\oplus Y
$$
The action of the group $\GL(V)$ on grassmannian on language
of angular operators is given by the formula
$$R\to (C+DR)(A+BR)^{-1}$$

  Let us consider a symmetric space $G/H$ and let us
realize it as a set $\cal S$. Call to mind that
a point of $\cal S$ is a ordered pair of subspaces $(Q,R)$.
Fix a pair $(X,Y)\in \cal S$.
For each point $(Q,R)\in \cal S$ we associate the pair of operators
 $$(M,N)$$
where

  $M:X\to Y$ is the angular operator of
the subspace $Q$
associated with the decomposition $X\oplus Y$ and

    $N:Y\to X$ is the angular operator of the subspace $R$
 associated with the decomposition $Y\oplus X$.

The condition $Q\cap R=0$ in our coordinates means
\begin{equation}
\det(1-MN)\ne0
\end{equation}

All other conditions also can be written in very simple form

\medskip

a) Consider the case then $Q$, $R$ (and in particular $X$, $Y$)
 are isotropic with respect to a form $B$ (Lists 1,2). 
For for all $x_1,x_2\in X $ we have $x_1+Mx_1,x_2+Mx_2\in Q$. Hence
\begin{multline*}
0= B(x_1+Mx_1, x_2+Mx_2)=\\
=B(x_1, x_2)+B(Mx_1, Mx_2) +B(Mx_1, x_2)+B(x_1, Mx_2)=\\
=0+0+B(Mx_1, x_2)+B(x_1, Mx_2)
\end{multline*}
(and similarly for $N$). Thus
\begin{equation}
B(Mx_1, x_2)+B(x_1, Mx_2)=0
\end{equation}
In matrix coordinates it means that a matrix is
symmetric, skewsymmetric, hermitian, antihermitian.
(depending on a type of the form $B$)

\medskip

b) Consider the case when $R$, $Q$ are orthogonal with
respect to a managing form $D$ (Lists 1,4). Then for all
$x\in Q$, $y\in R$
\begin{multline*}
0=D(x+Mx,y+Ny)=D(x,y)+D(Mx,Ny)+D(x,Ny)+D(Mx,y)=\\
=0+0+D(x,Ny)+D(Mx,y)
\end{multline*}
Thus
\begin{equation}
D(x,Ny)+D(Mx,y)=0
\end{equation}
In matrix language it give condition $N=\pm M^*$ or $N=\pm M^t$
(depending on the type of the form $D$).

\medskip

c) Consider the cases then  $R$ and $Q$ are linked by managing
semiinvolution
$J$(Lists 1,3). Then we obtain the condition
\begin{equation}
N=JMJ^{-1}
\end{equation}

\medskip

d) In the cases marked by $\star$ different symmetric spaces
$G/H_i$ are separated by the hypersurface
$\det(1-MN)=0$.

{\sc Remark.} Emphasis that equations (3.2)--(3.4)
are linear. For each point we $(X,Y)\in\cal S$
we constructed a map on the manifold $G/H=\cal S$. Thus we obtained atlas on the manifold $G/H$

{\sc Remark.} For riemann noncompact symmetric spaces our construction
is equivalent to realization of the type  "matrix ball"
(see for instance \cite{Ner}, Addendum A).

{\bf 3.5. Hua Loo Keng double ratio.} Let $(Q_1,Q_2)$, $(R_1,R_2)$ be
points of a symmetric space
${\cal S}=G/H$.
Let  $M:Q_1\to Q_2$ be the angular operator of the subspace
$R_1$ associated to the decomposition $V=Q_1\oplus Q_2$.
Let $N:Q_2\to Q_1$ be an angular operator of the subspace
$R_2$ associated to the decomposition $V=Q_2\oplus Q_1$.
Then $NM$ is a canonically defined operator $Q_1\to Q_1$.
Its eigenvalues $(\lambda_1,\lambda_1,\dots)$
are invariants of a pair of points
$(Q_1,Q_2)$, $(R_1,R_2)$ under the action of the group $G$.   This
construction is close to the usual double ratio of 4 points
of projective line. For several series of classical symmetric spaces
it was defined by Hua Loo Keng \cite{Hua}.

{\sc Remark.} Let $(X_1,X_2)$ be coordinates of a point $(Q_1,Q_2)$,
and  $(Y_1,Y_2)$ be  coordinates of a point$(R_1,R_2)$
in the sense of previous Subsection.  Then double ratio
coincides with eigenvalues
of the matrix
$$(1-Y_2X_1)^{-1}(X_2-Y_2)(1-Y_1 X_2)^{-1}(X_1-Y_1)$$

{\bf 3.6. Goncharov--Gindikin conformal structures.}
Le $\Gr^*$ be one of our grassmannians. Fix a point $P\in \Gr^*$
and an integer $\alpha=1,2,\dots,\dim P$.
By ${\rm D}_\alpha(P)$ we denote
the space of all subspaces $Q\in \Gr^*$
such that codimension of $P\cap Q$ in $P$ is less or equals $\alpha$.
By $T_P$ we denote the tangent space to $\Gr^*$ in a point
$P$.
Denote by ${\rm C}_\alpha(P)$ the cone  $T_P$ consisting of vectors
 tangent to variety ${\rm D}_\alpha(P)$. In this way for
each
$\alpha$ we obtaine a field of cones ${\rm C}_\alpha$
on grassmannian $\Gr^*$. Obviously this field of cones
is $G^*$-invariant.

 A.B.Goncharov and C.G.Gindikin (see \cite{Gin})
considered the field $C_1$
(or  $C_2$ if $C_1$ is empty).
It appeared that this structure
({\it in the case or rank $>1$})
"remember" the group $G^*$. precisely the pseudogroup
of diffeomorphisms of
$G/H$ preserving field of cones $C_1$ is the group
$G^*$ ( up to connected components).

\end{document}